\documentclass[11pt, reqno]{amsart}
\usepackage[latin1]{inputenc}
\usepackage{amssymb}
\usepackage{graphics}
\usepackage[all]{xy}

\usepackage{syntonly}
\usepackage{pstricks,subfigure,epsfig}
\usepackage{amsmath,amsfonts,amssymb,amsthm,euscript,upgreek}
\usepackage{enumerate}
\usepackage{multirow}
\usepackage{appendix}

\newtheorem{theoreme}{Theorem}[section]
\newtheorem{definition}[theoreme]{Definition}
\newtheorem{remarque}[theoreme]{Remark}
\newtheorem{proposition}[theoreme]{Proposition}
\newtheorem{corollaire}[theoreme]{Corollary}
\newtheorem{lemme}[theoreme]{Lemma}

\newtheorem{conjecture}{Conjecture}

\newtheorem{question}[theoreme]{Question}

\newcommand{\real}{\mathbb{R}}

\newcommand{\R}{\Bbb{R}}

\newcommand{\Vol}{\operatorname{Vol}}

\theoremstyle{plain}

\begin{document}

\title[Maximizing torsional rigidity on Riemannian manifolds] {Maximizing torsional rigidity on Riemannian manifolds}

\author{Lucio Cadeddu}

\address{Dipartimento di Matematica \\
         Universit\`a di Cagliari}
         \email{cadeddu@unica.it}
         
 \author{Sylvestre Gallot}
\address{Universit\'e Joseph Fourier Grenoble \\
          Institut Fourier}
\email{sylvestre.gallot@ujf-grenoble.fr}

\author{Andrea Loi}
\address{Dipartimento di Matematica \\
         Universit\`a di Cagliari}
         \email{loi@unica.it}

\thanks{
The second author was supported by the program  \lq\lq Visiting professor'' of the University of Cagliari, granted by the Province of Cagliari}
\thanks{
The third author was (partially) supported by ESF within the program \lq\lq Contact and Symplectic Topology''}

\date{\today}

\begin{abstract}
Let $\,(M,g)\,$ be a $n$-dimensional Riemannian manifold and $\,\Omega\,$ be any compact 
connected domain in $\,M$. We study the problem of finding the {\em maxima} of the functional
$\, {\mathcal E} (\Omega)\,$ (known as  {\em torsional rigidity} associated to $\Omega$)  among all domains of prescribed volume $v$. Our results show that for a given Riemannian manifold which is strictly isoperimetric at one of its points the maximum of such functional is realized by the geodesic ball centered at this point.
More generally, we prove estimates for the functional $\, {\mathcal E} (\Omega)\,$ by comparison with symmetrized domains. We also investigate on finding sharp upper bounds for the functional $\, {\mathcal E} (\Omega)\,$, under certain conditions on the geometry of $\,(M,g)\,$ and of $\Omega$. Finally we find an universal upper bound for $\, {\mathcal E} (\Omega)\,$ in terms of the isoperimetric Cheeger constant.
\end{abstract}

\maketitle

\section{Introduction}
Let $\,(M,g)\,$ be a $n$-dimensional Riemannian manifold (compact or not), $\,d\,$ be the associated Riemannian
distance and $\,dv_g\,$ the associated Riemannian measure. Let $\,\Omega\,$ be any compact 
connected domain in $\,M$, with {\it smooth boundary} $\,\partial \Omega\,$ (by this, in the case where
$\,M\,$ is compact, we also intend that the interior of $\,M \setminus \Omega\,$ is a non empty open set). 
Let us denote by $\,\Delta\,$ the Laplacian
\footnote{If  $f$ is a smooth function on $M$ then $\,\Delta \, f = - {\rm Trace} ( \nabla d f )\,$ and so  the Euclidean Laplacian writes
$\,\Delta = - \sum_{i = 1}^n \, \frac{\partial^2\, \ }{\partial^2 x_i}\,$.\\} on $\,M\,$ associated to  the 
Riemannian metric $\,g\,$, and let 
$\,f_\Omega\,$  be a  solution 
of the following Dirichlet problem

\begin{equation}\label{tempssortie}
\left\{\begin{array}{l}
\Delta \, f = 1 \ \ \ {\rm on}\ \ \  \Omega \ \ \\
\ \ \ f = 0 \  \ \  {\rm on}\    \  \partial \Omega .\ \ 
\end{array}\right. \,
\end{equation}
Let $\,C_c^\infty (\Omega) \,$ be the space of $\,C^\infty \,$ functions with compact 
support in the interior  of $\,\Omega\,$ and let $\,H_{1,c}^2 (\Omega)\,$ be its completion with respect to the Hilbert (Solobev) norm 
$\| f \|_{H_1^2(\Omega)} = (\| f \|^2_{L^2(\Omega)} + \| \nabla f \|^2_{L^2(\Omega)})^\frac{1}{2}$.  
As  $\,f_\Omega\,$  is regular

and vanishes on 
$\,\partial \Omega\,$,  it is  not hard to see that $\, f_\Omega\in H_{1,c}^2 (\Omega)\,$. Moreover, $f_\Omega(x) > 0$ for any $x\in \mathring\Omega$.

 On the space $\,H_{1,c}^2 (\Omega)\,$  let us consider the functional $\,E_\Omega\,$ 
defined by \footnote{See   Remark \ref{explanation}  in Section \ref{compsameman}  below for the explanation  of the appearance of  ${\rm Vol} (\Omega)$ in the definition of $E_\Omega$.}
\begin{equation}\label{EOmega}
E_\Omega (f) = \dfrac{1}{{\rm Vol} (\Omega)} \left(2 \int_\Omega \, f\ dv_g - 
\int_\Omega \,|\nabla f|^2 \ dv_g\right).
\end{equation}
Computing the first variation of $E_\Omega$ at the point $f$, we get
that 
$$(\, f\ {\rm is\ a\ critical\ point\ of}\ E_\Omega \,) \iff (\, f\ {\rm is\ a\ solution\ 
of\ (\ref{tempssortie})}\, ) .$$
The existence of (at least) one solution $\,f_\Omega\,$ of (\ref{tempssortie}) proves 
that the functional  $\,f \mapsto E_\Omega (f)\,$ (defined on $\,H_{1,c}^2 (\Omega)\,$) admits at least one critical point.
Moreover, the functional $\,E_\Omega\,$ is  strictly concave 
hence it admits a unique critical point, which is its (unique) absolute maximum and, consequently,  problem (\ref{tempssortie}) admits $\,f_\Omega\,$ as a unique solution.

\vskip 0.3cm

Therefore one can give  the following  definition  (see \cite{MD} and references therein for more details).

\begin{definition}\label{fonctionnelle}
Let $\Omega\subset M$ as above. 
The {\em torsional rigidity}
 \footnote{When
$\,\Omega\,$ is a domain of the Euclidean plane, 
${\mathcal E} (\Omega) \,$ is the torsional rigidity of a beam whose cross section 
is $\,\Omega\,$.\\}  of $\Omega$ is the value
$${\mathcal E} (\Omega) = E_\Omega\, (f_\Omega) = \max_{f \in H_{1,c}^2 (\Omega)}
\left( E_\Omega\, (f)\right) .$$
\end{definition}

Consider  the functional ${\mathcal E}:\Omega\rightarrow \R$ restricted to the set of domains
$\Omega\subset M$ with smooth boundary and prescribed volume $v$. It is known that
its critical points are  the {\em harmonic domains} 
namely those domains $\Omega\subset M$ such that 
the function $\|\nabla f_\Omega (x)\|$
 is constant on the boundary $\,\partial \Omega$.
In the literature, the proof of this assertion is based on a Brownian motion probabilistic
argument (see \cite{MD}, Proposition 2.1).
Recall also   the classical fact  that,  in a 
Riemannian manifold $(M,g)$  which is harmonic at one of its point $x_0$  (see next section for the definition), 
every geodesic ball  centered at $x_0$ is  a harmonic domain. Hence one has  the following fundamental question.

\begin{question}\label{question2}
On  a Riemannian  manifold $(M, g)$ which is harmonic at some of its  points $\,x_0\,$
is every harmonic domain a geodesic ball centered at $\,x_0\,$?
\end{question}

With respect to this question, when $(M, g)$ is a real space form (and hence harmonic at each of its points) the following results hold true.
J. Serrin (\cite{Se}) proved that every harmonic domain 
of $\,(\mathbb R^n , {\rm can}.)\,$ is a ball. S. Kumaresan and J. Prajapat (\cite{K-P}) extended Serrin's result, proving that every harmonic domain
of $\,(\mathbb R \rm H^n, can.)\,$ is a geodesic ball and that every harmonic domain
of $\,(\mathbb S^n, can.)\,$ whose closure is contained in an hemisphere is a geodesic ball.
On the other hand  it is not true that every harmonic domain in $\,(\mathbb S^n, can.)\,$ is a geodesic ball
(semi-classical counterexamples are given by tubular neighbourhoods, in $S^3$, of some geodesic circle $S^1$ and, more generally, by domains with isoparamteric boundary in $\mathbb S^n$) and so the previous question has, in general, a negative answer.

\vskip 0.3cm
One of the aims of  the present  paper is to address the problem of studying the {\em maxima}  (instead of the critical points) of the functional
$\,\Omega \mapsto {\mathcal E} (\Omega)\,$ 
among all domains of prescribed volume $v$ (obviously such maxima are harmonic domains).

The following  represents our first result.
It shows that, for a given Riemannian manifold $(M, g)$ which is  {\em strictly  isoperimetric} at one of its points (see Definition \ref{isoperimetre}
 in Section \ref{harmonic} below), the geodesic ball centered at this point  realizes the  maximum  of the torsional rigidity ${\mathcal E} (\Omega)$ 
and moreover a maximum is a geodesic ball.

\begin{theoreme}\label{comparaison1}
Let $(M, g)$ be a Riemannian manifold  which is  isoperimetric at some point $x_0 \in M$ and let $\Omega$ be any compact domain with smooth boundary in $M$.
Let  $\Omega^*$  be  the geodesic ball of $\,(M, g)$ centered at $x_0$ such that ${\rm Vol} (\Omega^*) = {\rm Vol} (\Omega)$ ,
then $${\mathcal E} (\Omega) \le {\mathcal E} (\Omega^*).$$
Moreover, if $(M, g)$ is strictly isoperimetric at $x_0$ then the equality 
$\,{\mathcal E} (\Omega) = {\mathcal E} (\Omega^*)\,$ is realized if and only if $\,\Omega\,$ is isometric to $\,\Omega^*\,$.
\end{theoreme}

An immediate consequence of Theorem \ref{comparaison1} is that, on the Euclidean space, on the Hyperbolic space and on the 
canonical sphere, the geodesic balls are the domains which realize the maximum
of the functional $\,\Omega \mapsto {\mathcal E} (\Omega)\,$ among all domains of prescribed 
volume $\,v\,$ (this was already known for domains in  the Euclidean space, in the Hyperbolic space 
and in the canonical open hemisphere, see \cite{MD} and also \cite{Ca}). The present result extends this property
to the whole sphere and to domains in some manifolds of revolution described in Section \ref{Example}
below.

The proof of Theorem \ref{comparaison1} is a particular case of the following Theorem \ref{comparaison}. 
The reader is referred to Section \ref{compsameman} below for the definition (Definition \ref{espacemodele}) 
of {\em symmetrized domain}  $\Omega^*$ of a given domain $\Omega$  and  of {\em pointed isoperimetric model space} $(M^*, g^*, x^*)$
 of a Riemannian manifold $(M, g)$ (PIMS in the sequel).

\begin{theoreme}\label{comparaison}
Let $\,(M,g)\,$ be a Riemannian manifold and let $\,(M^*,g^*, x^*)\,$ be a PIMS
for $\,(M,g)\,$.
Let $\,\Omega \,$ be any compact domain with smooth boundary in $\, M \,$ 
and  let $\,\Omega^*\,$ 
be its symmetrized domain.
Then
$${\mathcal E} (\Omega) \le {\mathcal E} (\Omega^*).$$
Moreover, if $\,(M^*,g^*, x^*)\,$ is a strict PIMS
for $\,(M,g)\,$ then the equality 
$\,{\mathcal E} (\Omega) = {\mathcal E} (\Omega^*)\,$ is realized if and only if $\,\Omega\,$
is isometric to $\,\Omega^*\,$.
\end{theoreme}

From Theorem \ref{comparaison} the following question naturally arises.

\begin{question}\label{question3}
Can we find sharp universal upper bounds C(v) for 
the torsional rigidity $\,{\mathcal E} (\Omega) \,$  which are
independent on the geometry of $\,(M,g)\,$ (except for some a priori bounds on its 
curvature and diameter) and on the geometry of the domain $\,\Omega \subset M\,$ (provided
that this domain has prescribed volume $v$)?
\end{question}

In order to attack this question one needs to find a unique \lq\lq universal'  strict PIMS for all the Riemannian manifolds that belong to a given class.

\vskip 0.3cm

In the noncompact case one has the following  well-known conjecture which is called Cartan--Hadamard's conjecture (or Aubin's conjecture) in the literature. We recall that a 
 Cartan--Hadamard manifold is  a complete  simply connected Riemannian manifold with non positive sectional curvature.

\begin{conjecture}\label{Croke}
The Euclidean n-dimensional space $E^n$, pointed in any point $x^* \in E^n$, is a
strict PIMS  for every Cartan-Hadamard manifold.
\end{conjecture}

This conjecture is known to be true when the dimension $\,n \,$ is equal to $\,2\,$
(it is a classical fact, using the Gauss-Bonnet formula, proved for the first time by A. Weil in \cite{Weil}),
in dimension $\,4\,$ (it was proved by C. B. Croke \cite{Cr}, using Santalo's formula) and in dimension
$\,3\,$ (it is a more recent proof by B. Kleiner \cite{Kl}).  In higher dimensions, Conjecture \ref{Croke} is still open.

Using  these results  
we immediately get the following corollary
of  Theorem \ref{comparaison} which provides an answer to Question \ref{question3} in the noncompact case.

\begin{corollaire}\label{comparaisonCroke}
Let $(M, g)$ be a  Cartan--Hadamard manifold of  dimension $n\le 4$.
For every compact domain $\Omega\subset M$ with smooth boundary, one has 
$${\mathcal E} (\Omega) \le {\mathcal E} (\Omega^*),$$
where $\,\Omega^*\,$ is the Euclidean ball with the same volume as $\,\Omega\,$.
Moreover the equality $\,{\mathcal E} (\Omega) = {\mathcal E} (\Omega^*)\,$ is realized if and only if 
$\,\Omega\,$ is isometric to an Euclidean ball. 
\end{corollaire}

Notice that, if the conjecture \ref{Croke} was true in every dimension $n$, the Corollary \ref{comparaisonCroke} would be automatically true in any dimension.

In the compact  case one has the celebrated Gromov's isoperimetric inequality and 
its generalization due to  P. B\'erard, G. Besson and S. Gallot  (see respectively Theorems \ref{Gromov} and  \ref{B-C-G}).
Using these isoperimetric results  and a Theorem of G. Perelman (Theorem \ref{Perelman})  we obtain the following result 
which gives an answer to Question \ref{question3} in the compact case.
Moreover it shows  that a manifold has the same geometry and topology as a sphere by the knowledge of the value of the torsional
rigidity of one of its domains.

\begin{theoreme}\label{compsphere2}
For every complete, connected Riemannian manifold $\,(M,g)\,$ whose Ricci
curvature satisfies $\,{\rm Ric}_g \ge (n-1).g$, for every compact domain with smooth boundary 
$\,\Omega \,$ in $\, M \,$, let $\,\Omega^*\,$ be a geodesic ball of  the canonical sphere 
$\, (\mathbb S^n , g_0 )\,$ such that 
$\,\dfrac{{\rm Vol} (\Omega^*, g_0)}{{\rm Vol} (\mathbb S^n , g_0 )} 
= \dfrac{{\rm Vol} (\Omega , g )}{{\rm Vol} (M , g )}$, 
then 
$${\mathcal E} (\Omega) \le {\mathcal E} (\Omega^*).$$
Morever, 
\begin{itemize}
\item [(i)] 
if there exists some domain $\,\Omega \subset M\,$ such that $\,{\mathcal E} (\Omega) = 
{\mathcal E} (\Omega^*)\,$  then $\, (M, g) \,$ is isometric to $\,(\mathbb S^n , g_0 )\,$ and
$\,\Omega\,$ is isometric to $\,\Omega^*\,$.
\item [(ii)] If there exists some domain $\,\Omega \subset M\,$ such that 
$$\,{\mathcal E} (\Omega) >
 \left( 1 - \dfrac{\int_0^{\frac{\varepsilon (n , \kappa)}{2}}\, (\sin t)^{n-1}\, dt}
{\int_0^{\frac{\pi}{2}}\, (\sin t)^{n-1}\, dt}\right)^{\frac{2}{n}}{\mathcal E} (\Omega^*),$$
(where $\,-\kappa^2\,$ is a lower bound for the sectional curvature of $\,(M,g)\,$ and where
$\,\varepsilon (n , \kappa)\,$ is the Perelman constant  described in Theorem \ref{Perelman})
then $\, M \,$ is diffeomorphic to $\,\mathbb S^n \,$.
\end{itemize}
\end{theoreme}

Our last result is the following  Theorem \ref{cheeger}, where we provide  a sharp universal bound for the torsional rigidity of any domain of a compact Riemannian manifold  $(M, g)$ in terms of its Cheeger  isoperimetric constant  $H(M, g)$.  Let us recall that  $H(M,g)$  is defined by 
$$H(M,g) = \inf_\Omega \left( \dfrac{{\rm Vol}_{n-1} (\partial \Omega)}{\min \left({\rm Vol} (\Omega),
{\rm Vol} (M \setminus \Omega)\right)}\right),$$
where $\,\Omega\,$ runs in the set of all domains with smooth boundary in $M$. 

\begin{theoreme}\label{cheeger}
Let $(M,g)$ be any compact Riemannian manifold and  let $\Omega$ be any compact 
domain with smooth boundary in $M$ such that  $Vol (\Omega) \le \frac{1}{2} {\rm \ Vol} (M)$.
Then ${\mathcal E} (\Omega) \le  \dfrac{1}{H(M,g)^2} .$\\
\end{theoreme}

This inequality is sharp:  at the end of the paper we shall exhibit examples of sequences of  Riemannian manifolds $(M, g_\epsilon)$, $0<\epsilon <1$ 
and of domains $\Omega_\epsilon \subset M$ such that ${\mathcal E} (\Omega_\epsilon)H^2(M , g_\epsilon) \rightarrow 1$  as $\epsilon \rightarrow 0.$ 

\vskip 0.3cm

The paper is organized as follows.
In Section 2 we recall the definition of manifolds which are harmonic and  isoperimetric at a point and provide examples of non standard isoperimetric Riemannian manifolds.
 In Section 3 we give the definition of PIMS for a given manifold $(M,g)$ and we prove Theorems 1.3 and 1.4. 
The main tool in the proof of this theorem is the Theorem of symmetrization (Theorem \ref{symetrisation}) which 
 gives precise relationships between the integrals that appear in the definition of ${\mathcal E} (\Omega)$ when calculated on $\Omega$ and on its symmetrized $\Omega^*$. 
In Section 4 we investigate how to compare torsional rigidities of domains of two different compact manifolds. This will allow us to prove Theorem 1.7. This section  ends with the proof of Theorem \ref{cheeger} and of its sharpness  (Proposition \ref{propsharp}).

\section{Harmonic and isoperimetric manifolds at one point}
We briefly recall the definition of harmonic manifolds (the reader is referred to \cite{Be}, \cite{R-D} and references therein for details).
Let us recall that there exists, in any Riemannian manifold $\, (M,g)\,$, and for any point
$\,x_0 \in M\,$, a closed subset of measure zero, the {\em cut-locus}  of $\,x_0\,$ (denoted
by $\,{\rm Cut} (x_0)\,$) such that
the exponential map $\,{\rm exp}_{x_0}\,$ is a diffeomorphism from an open subset
$\, U_{x_0}\,$ of the tangent space $\,T_{x_0} M\,$ onto $\,M \setminus {\rm Cut} (x_0)\,$. 
Let $\,\mathbb S_{x_0}\,$ be the (Euclidean) unit sphere of the Euclidean space 
$\,\left(T_{x_0} M , g_{x_0}\right) \,$ and let us define the open subset 
$\,\widetilde U_{x_0} \subset ]\,0 \,,\,+ \infty\,[ \times \mathbb S_{x_0}\,$ 
as the pull-back of $\, U_{x_0}\,$ by the map $\, (t, v) \mapsto t.v\,$ from 
$\,]\,0 \,,\,+ \infty\,[ \times \mathbb S_{x_0}\,$ to $\,T_{x_0} M\,$; this provides a 
generalization of the usual \lq \lq polar coordinates" by the notion of  {\em normal coordinates} 
$$\phi : \left\{\begin{array}{l}
\widetilde U_{x_0} \to  U_{x_0} \to M \setminus {\rm Cut} (x_0)\\
(t, v) \mapsto t.v \mapsto {\rm exp}_{x_0} (t.v) \end{array}\right.$$
In this coordinates system, let us write the Riemannian measure at the point $\,\phi (t,v)\,$ as
\begin{equation}\label{theta}
\phi^* dv_g = \theta (t,v)\, dt\,dv,
\end{equation}
where $\,dv\,$ is the canonical measure of the canonical sphere $\,\mathbb S_{x_0}\,$.
This defines  $\,\theta (t,v)\,$ as the density of the measure $\, \phi^* dv_g \,$ with respect to
the measure $\, dt\,dv \,$.

\begin{definition}\label{harmonicpointe} 
Let $\, (M, g)\,$ be a Riemannian manifold and $\,x_0 \,$ be a point in $\,M\,$,
$(M,g)$, is said to be {\em harmonic at $\,x_0 \,$}  iff the following two conditions
are satisfied:
\begin{itemize}
\item $\, U_{x_0}$ is equal  to $T_{x_0}M$  or to an open   ball of the Euclidean space $\,\left(T_{x_0} M , g_{x_0}\right) \,$
{\rm (}and thus there exists some $\,\beta \in ] 0 , + \infty]\,$ such that 
$\,\widetilde U_{x_0} = ]\,0 \,,\,\beta\,[ \times \mathbb S_{x_0}\,${\rm )}. 

\item for every $\,t \in ]\,0 \,,\,\beta\,[\,$, $\,\theta (t,v)\,$ does not depend on $\,v\,$ and will be, in this case, denoted by $\theta(t)$.
\end{itemize}
\end{definition}

\begin{definition}\label{harmonic}
A Riemannian manifold $\, (M,g)\,$ is said to be {\em harmonic} iff it is
harmonic at each of its points.
\end{definition}

For example, spaces of revolution are  harmonic at their pole(s),
but they are generally not harmonic in the sense of the Definition \ref{harmonic}. See Section  \ref{Example} below.

\begin{definition}\label{isoperimetre}
Let $\,(M, g)\,$ be a Riemannian manifold and $\,x_0\,$ a point of $M$. 
The manifold $\,(M,g)$ is said to be {\em isoperimetric at $x_0$}  if it is  harmonic at $x_0$  and if, for any compact
domain $\,\Omega \subset M\,$ with smooth boundary, the geodesic ball 
\footnote{Such a domain $\Omega^*$ writes $B(x_0 , R_0)$, where $R_0$ is
the solution of the equation ${\rm Vol} \left(B(x_0 , r ) \right) =  {\rm Vol} \left( \Omega \right)$, this solution 
exists and is unique because, being $(M, g)$ harmonic at $x_0$, the function 
$r\mapsto {\rm Vol} \left(B(x_0 , r ) \right)={\rm Vol}_{n-1}(S^{n-1})\int_0^r\theta (t)dt$
 is a continuous strictly increasing function. We do not really need harmonicity to prove that but, in this case, the proof is simpler.} 
$\Omega^*$ centered at $x_0$ with the same volume as $\Omega$ satisfies $\, {\rm Vol}_{n-1}(\partial \Omega^*) \le {\rm Vol}_{n-1}(\partial \Omega)\,$;
the same manifold is said to be {\em strictly isoperimetric at $\,x_0\,$}  if, moreover, the equality occurs iff $\,\Omega\,$ is isometric to $\Omega^*$.
\end{definition}

The Euclidean space, the  Hyperbolic Space and the Sphere are strictly isoperimetric 
at every point  (for proofs\footnote{A short history of these
proofs is given in \cite{B-Z}, Section 10.4, see also Section 8.6.\\} see for instance \cite{B-Z}).
These examples  are the only known examples (up to homotheties) of Riemannian manifolds 
which are isoperimetric at every point.
If we only require the  Riemannian manifolds to be isoperimetric at (at least)
one point, we get much more examples. 
In fact, some non standard spaces of revolution 
are isoperimetric at one pole, as in the following  example.

\subsection{Examples of nonstandard Riemannian manifolds which are isoperimetric at some point} \label{Example}
A (noncompact) space of revolution $\,(M,g)\,$ with
only one pole $x_0$ is
such that $ \left(M \setminus \{x_0\} , g\right)$ is isometric to 
$\,]\,0\,, \, + \infty \,[ \times \mathbb S^{n-1}$, endowed with a Riemannian metric
of the type $\, (dt)^2 + \,b(t)^2\, g_{\mathbb S^{n-1}}\,$, where $\,b\,$ is a smooth
strictly positive function whose extension to $\,[0\,, \, + \infty \,[ \,$ satisfies
$\,b(0) = 0\,$ (and $\,b'(0) = 1\,$ if we want the metric
to be regular at $\,x_0 $), where $\, g_{\mathbb S^{n-1}}\,$ is the canonical metric of the sphere
$\,\mathbb S^{n-1}\,$, and where $\, \{ 0\} \times \mathbb S^{n-1}\,$ is identified with the point
$\,x_0\,$.\\
A (compact) space of revolution $\,(M,g)\,$ with two poles $\,x_0\,$ and $\,x_1\,$ is
such that $\, \left(M \setminus \{x_0 , x_1\} , g\right)\,$ is isometric to 
$\,]\,0\,, \, L \,[ \times \mathbb S^{n-1}$, endowed with a Riemannian metric
of the type $\, (dt)^2 + \,b(t)^2\, g_{\mathbb S^{n-1}}\,$, where $\,b\,$ is a smooth
strictly positive function whose extension to $\, [0\,, \, L ]\,$ satisfies
$\,b(0) = b(L) = 0\,$ (and $\,b'(0) = 1\,$,
$\,b'(L) = -1\,$ if we want the metric to be regular at $\,x_0\,$ and $\,x_1\,$) and where 
$\, \{ 0\} \times \mathbb S^{n-1}\,$ (resp. $\, \{ L\} \times \mathbb S^{n-1}\,$) is 
identified with the point $\,x_0\,$ (resp.  with the point $x_1$). It is not hard to see that a  space of revolution is  harmonic at 
any of its poles.

We now show that some (nonstandard)  spaces of revolution are isoperimetrtic at their poles. 
The first  example is given by a 2-dimensional cylinder $\,[0 , +\infty \,[ \times \mathbb S^1\,$ 
(resp. $\,[0 , L] \times \mathbb S^1\,$) with 1 hemisphere glued to the boundary 
$\,\{0\} \times \mathbb S^1\,$ (resp. with 2 hemispheres respectively glued to the boundaries
$\,\{0\} \times \mathbb S^1\,$ and $\,\{L\} \times \mathbb S^1\,$). 
Other examples are given by  the paraboloid of revolution $z=x^2+y^2$
or the hyperboloid of equation $x^2 + y^2 - z^2 = -1,\  z>0$ in 
$\,\mathbb R^3$ (isoperimetric at their pole).

More generally, a large classe of nonstandard examples is given by the 

\begin{theoreme}
(\cite{MO}, Theorem 1.2) Consider the plane $\R^2$ equipped with a complete and   rotationally invariant Riemannian metric $g$ such that the Gauss curvature is positive and a strictly decreasing function of the distance  from the origin. Then $(\R^2, g)$ is isoperimetric at  the origin.
\end{theoreme}
\begin{remarque}\rm
Notice  that it is not true that every space of revolution 
is isoperimetric at its pole: for example 
let us consider the hypersurface of revolution $S$ in $\R^3$ of equation $\, x^2 + y^2 + ( |z| + \cos R)^2 = 1\,$, whose poles are $x_0 =(0,0,1-\cos R)$ and $x_1= -x_0$; then the geodesic ball $\,B(x_0 , R) \,$ is the subset $\, \{(x,y,z) \in S:\ z>0\}\,$ and $\,\partial B(x_0 , R) \,$
is the circle $\, x^2 + y^2 = \sin^2 R \ ,\ z=0\,$ whose length is $\,2\pi \sin (R)\,$.
The plane $\,y = 0\,$ separates the surface $\,S\,$ in two symmetric domains 
 $\,\Omega_1\,$ and $\,\Omega_2\,$, which have the same area as $\,B(x_0 , R) \,$ and whose boundary is the union of two arcs of circle of length $2 R$, we thus have 
$${\rm length} (\partial \Omega_1) = 4 R < 2\pi \sin (R) = {\rm length} (\partial B(x_0 , R)).$$
\end{remarque}

\section{The symmetrization of a function, the  theorem of symmetrization and the  proofs of Theorems \ref{comparaison1} and \ref{comparaison}}\label{compsameman}
Let $\,(M,g)\,$ and $\,(M^*,g^*)\,$ be two Riemannian manifolds such that $\,{\rm Vol}(M,g)\,$ 
and $\,{\rm Vol}(M^*,g^*)\,$ are both infinite or both finite. Let us define the constant 
$\,\alpha (M, M^*)\,$ by

$$\alpha (M, M^*) =  \left\{\begin{array}{l}
\ \ \ \ \ 1 \ \ \ \ \ \ \ \ \ \ \  {\rm if \  Vol} (M,g) \  {\rm and \  Vol} (M^*,g^*) \ 
{\rm are\ both\ infinite},\\  \\ 
\dfrac{{\rm Vol}(M,g)}{{\rm Vol}(M^*,g^*)} \ \ \ \  {\rm if \  Vol} (M,g) \  
{\rm and \  Vol} (M^*,g^*) \ {\rm are\ both\ finite.}\end{array}\right.$$

\begin{definition}\label{espacemodele}
Let $x^*$ be a fixed point of $M^*$.

\noindent
a) For any compact domain $\,\Omega \subset M\,$ with smooth boundary, one defines
its {\em symmetrized domain}
\footnote{The symmetrized domain $\Omega^* = B(x^* , R)$, where $R$ is
the solution of the equation $\, {\rm Vol} \left(B(x^* , r ) \right) = \alpha (M,M^*)^{-1}\,  
{\rm Vol} \left( \Omega \right)$, always exists and is unique because  
$\,\alpha (M,M^*)^{-1}\,  {\rm Vol} \left( \Omega \right)  \in\  ]\,0, {\rm Vol} (M^*,g^*)\,[ \,$
and $\, r \mapsto {\rm Vol} \left(B(x^* , r ) \right) \,$ is a continuous
strictly increasing function whose image is the closure of $\,]\,0, {\rm Vol} (M^*,g^*)\,[\,$ in $\,]\,0, +\infty \,[\,$.\\}
$\Omega^*$ (around the point $x^*$) 
as the geodesic ball  of $\,(M^*,g^*)\,$,
centered at $x^*$, such that $\, {\rm Vol} (\Omega^*) = \alpha (M,M^*)^{-1}\, 
{\rm Vol} (\Omega)\,$.

\noindent
b) $\,(M^*,g^*, x^*)\,$ is said to be a {\em pointed isoperimetric  model space}  ({\em PIMS) 
for $\,(M,g)\,$} if, 
for any compact domain $\,\Omega \subset M$, with smooth boundary, the symmetrized domain $\Omega^*$ satisfies the {\em isoperimetric inequality} 
$\, {\rm Vol}_{n-1}(\partial \Omega) \ge  \alpha (M,M^*)\, {\rm Vol}_{n-1}(\partial \Omega^*)$;
the same manifold is said to be a {\em strict PIMS}  if, moreover, the equality occurs iff $\,\Omega\,$ is {\em isometric} to $(\Omega^*, g^*)$.
\end{definition}

\begin{remarque}\rm\label{Volrelatifs} When the two manifolds have different finite volumes
(i.e. when $\,\alpha (M, M^*) \neq 1\,$), we are compelled
to make the assumption $\, {\rm Vol} (\Omega^*) = \alpha (M,M^*)^{-1}\, {\rm Vol} (\Omega)\,$
(which, in this case means that the relative volumes 
${\rm Vol}(\Omega)/ {\rm Vol}(M,g)$  and ${\rm Vol}(\Omega^*)/ {\rm Vol}(M^*,g^*)$ are equal) instead of the usual
assumption $ {\rm Vol} (\Omega^*) = {\rm Vol} (\Omega)\,$.
In fact, if the symmetrized domain $\Omega^*$ is defined by the equality ${\rm Vol} (\Omega^{*}) = {\rm Vol} (\Omega)$,
it is hopeless to expect some bound from below for $\, {\rm Vol}_{n-1}(\partial \Omega)\,$ in terms of 
$\, {\rm Vol}_{n-1}(\partial \Omega^{*})\,$. Indeed, if 
$\, {\rm Vol} (M, g) > {\rm Vol} (M^*,g^*)\,$,
then $\,\Omega^{*}\,$ would not exist when $\,\Omega\,$ is such that
$\, {\rm Vol} (M^*,g^*) < {\rm Vol} (\Omega) < {\rm Vol} (M,g)\,$; on the other hand, 
if $\, {\rm Vol} (M, g) < {\rm Vol} (M^*,g^*)\,$ and if 
$\,\Omega = M \setminus B(x , \varepsilon)\,$, where $\varepsilon\,$ is arbitrarily small and
$\,B(x , \varepsilon)\,$ is any geodesic ball of radius $\varepsilon\,$ in $\,(M,g)\,$, we get that
 ${\rm Vol}(M^* \setminus \Omega^*) = {\rm Vol}(M^*, g^*) - {\rm Vol}(M, g) + O(\epsilon^n)$ and ${\rm Vol}(\partial\Omega^*)$ does not go to zero 
while ${\rm Vol}(\partial \Omega) = O(\epsilon^{n-1})$ goes to zero when $\epsilon \rightarrow 0$. It is thus impossible to bound ${\rm Vol}_{n-1}(\partial\Omega)$ from below in terms of ${\rm Vol}_{n-1}(\partial\Omega^*)$
\end{remarque}

Let $\,(M,g)\,$ be a Riemannian manifold  and $(M^*,g^*, x^*)$ be a PIMS for $(M, g)$.
 Let $\,\Omega\subset M$ be any compact
domain with smooth boundary.  Let $\,\Omega^*\,$ be its {\em symmetrized domain}
in the sense of Definition \ref{espacemodele}.

Let $\,f\,$ be any smooth nonnegative function
on $\,\Omega \,$ which vanishes on the boundary, we denote by $\,\Omega_t\,$ (or by
$\,\{f>t\}\,$) the set of points $\,x \in \Omega\,$ such that $\,f(x) > t\,$. Let us
denote by $\,\{f = t\}\,$ the set of points $\,x \in \Omega\,$ such that $\,f(x) = t\,$;
notice that the set of critical points of $\,f\,$ is compact and thus its image $\,{\mathcal S} (f)\,$
by $\,f\,$ is compact and, by Sard's theorem, it has Lebesgue measure zero in
$\, [0 , \sup f]\,$. For any {\em  regular value} $\,t\,$ of $\,f\,$, namely for any   $t\in  [0 , \sup f] \setminus {\mathcal S} (f)\,$ the set 
$\,\{f = t\}\,$ is a smooth submanifold of codimension $\,1\,$ in $M$, which is
equal to $\,\partial \Omega_t\,$.
For any $\, t \in [0 , \sup f [\,$, let us define $\,\Omega_t^*\,$ as the symmetrized 
domain of $\,\Omega_t\,$,  i.e. the geodesic
ball $\,B\left(x^* , R(t)\right) \,$ whose radius $\,R(t) \,$ is chosen in such a way that
$\, {\rm Vol} \left(B\left(x^* , R(t)\right) \right) = \alpha (M,M^*)^{-1}\,{\rm Vol} (\Omega_t)\,$.
When $\,t = \sup f \,$, then $\,\Omega_{\sup f }\,$ is empty, and thus $\,R(\sup f ) = 0\,$.

The function $\, t \mapsto A(t) := {\rm Vol} (\Omega_t)\,$ is strictly decreasing because,
when $\,0 \le t < t' \le \sup f\,$, the set $\,\{ x \in X \ :\ t <f(x) < t'\}\,$ is a nonempty open set of 
nonzero volume; a consequence is that the function $\, t \mapsto R(t)\,$ is also strictly
decreasing
\footnote{Thus $\, t \mapsto R(t)\,$ is well defined (for every $\,t\,$) and injective.
However, it is generally not surjective nor continuous, moreover the measure of the set
$\,[0 , R_0] \setminus {\rm Image} (R)\,$ is generally not zero. This is one of the main problems when studying the regularity of $\,\bar f\,$,
and thus of $\,f^*\,$.\\}.

We then define $\, f^* : \Omega^* \to \mathbb R^+\,$, the \emph{symmetrized of $f$} in such a way that $\{\,f^* > t\}\, = \Omega^*_t$, namely, 
we decide that $\,f^* = \bar f \circ \rho\,$, where
$\,\rho (x) = d^* (x^* , x)\,$, where $\,d^*\,$ is the Riemannian distance on $\,M^*\,$ associated
to $\,g^*\,$, and where $\,\bar f : [0 , R_0] \to [0 , \sup f]\,$ is defined by
$$\bar f (r) := \inf \left( R^{-1} ([0 , r ]) \right) =\inf \{t\in [0 , \sup f]\, :\, R(t) \le r \}$$
$$ = 
\inf \{t\, :\, A(t) \le \alpha (M,M^*) \,{\rm Vol} B(x^* , r) \}$$

We now state the Theorem of symmetrization which represents the main tool
for the proof of our main results. Symmetrization methods have their origin in J. Steiner's works. The following classical application to functional analysis (also called rearrangements) generalizes to Riemannian manifolds ideas of G. Talenti.

\begin{theoreme}\label{symetrisation}
Let $(M,g)$ be a Riemannian manifold and  $(M^*,g^*, x^*)$ be a PIMS
for $(M,g)$. Let  $\Omega\subset M$ be  a compact domain with smooth boundary and 
$f$ be any smooth nonnegative function on  $\Omega$ which vanishes on its boundary.
 Let $f^*$ be the symmetrized function, constructed as above on the symmetrized geodesic ball $\,\Omega^* \,$ of  $\,(M^*,g^*)\,$, centered at the point
 $x^*$. Then
\begin{itemize}
\item[(i)] $\,f^*\,$ is Lipschitz (with Lipschitz constant $\,\|\nabla f\|_{L^\infty}$) and thus $\,f^*\,$ 
lies in $\, H_{1,c}^2\left( \Omega^* , g^*\right)\,$,
\item[(ii)] $\,\dfrac{1}{{\rm Vol} (\Omega)} 
\int_\Omega f(x)^p  \,dv_g (x) = \dfrac{1}{{\rm Vol} (\Omega^*)} 
\int_{\Omega^*} \left(f^* (x) \right)^p \, dv_{g^*} (x)\,$  for every $\,p \in [ 1 , +\infty [\,$,
\item[(iii)] $\, \dfrac{1}{{\rm Vol} (\Omega)}  \int_\Omega \|\nabla f (x)\|^2 \,dv_g (x) \ge
\dfrac{1}{{\rm Vol} (\Omega^*)} \int_{\Omega^*} \|\nabla f^* (x)\|^2 \,dv_{g^*} (x) \,$. If, moreover, $(M^*,g^*)$ is a strict PIMS for $(M,g)$ then equality holds iff the set $\{f>0\} \subset (\Omega, g)$ is isometric to the set $\{f^*>0\} \subset (\Omega^*, g^*)$
\end{itemize}
\end{theoreme}
The proof of this theorem can be obtained by following the same lines as in \cite{Ber} and \cite{Ga2} (see also \cite{bandle}).
We point out that one of the main tools in the proof of the theorem of symmetrization given in  \cite{Ber} and \cite{Ga2}, namely the {\em coarea formula}, is not entirely correct in these references, though this has no consequences on the results proved in these papers.  The  correct version of the coarea formula (see for instance \cite{B-Z}), pp. 104-107) is the following:

\vskip 0.3cm

\noindent
{\bf Coarea formula:} 
{\em Let $\,(M,g)\,$ be a Riemannian manifold 
and let $f : M \to \mathbb R$ be a smooth function.
Then, for any measurable function $\varphi$ on $M$, one has\footnote {By $\, \int_{\inf f}^{\sup f} \,$, we intend the integral 
on $]\inf f ,\sup f [ \setminus  {\mathcal S}(f)$ (because ${\mathcal S}(f)$ has measure zero).
 Moreover, as we only integrate with 
respect to regular values $\,t\,$ of $\,f\,$, $\,\{ f = t\}\,$ is a submanifold of 
codimension $1$ in $M$ and $da_t $ is  the $(n-1)$-dimensional 
Riemannian measure on $\{ f = t\}$ (viewed as a Riemannian submanifold of $(M , g)$).}
\begin{equation}\label{coarea}
\int_M \varphi (x)\ \|\nabla f (x)\|  dv_g (x) = \int_{\inf f}^{\sup f} \left( \int_{f^{-1}(\{t\})} 
\varphi (x) da_t (x)\right) dt.
\end{equation}}

Formula (\ref{coarea}) together with the corrected proof of the theorem of symmetrization will appear in a forthcoming survey paper.

\vskip 0.3cm

\noindent {\it Proof of Theorem \ref{comparaison}}\,:
\noindent
Let $\,f_\Omega\,$ be the unique solution of the problem (\ref{tempssortie}) on the domain 
$\,\Omega$, let $\,(f_\Omega)^*\,$ be the corresponding symmetrized function. By Theorem of Symmetrization \ref{symetrisation} (ii) and (iii) we get
$${\mathcal E} (\Omega) = E_\Omega ( f_\Omega) = \dfrac{1}{{\rm Vol} (\Omega)} 
\left( 2 \int_\Omega \, f_\Omega\ dv_g - \int_\Omega \,|\nabla f_\Omega|^2 \ dv_g \right)$$
$$ \le \dfrac{1}{{\rm Vol} (\Omega^*)} 
\left( 2 \int_{\Omega^*} \, (f_\Omega)^*\ dv_{g^*} - \int_{\Omega^*} \,|\nabla (f_\Omega)^*|^2 
\ dv_{g^*} \right) = E_{\Omega^*} \left( (f_\Omega)^*\right) .$$
Let us recall that the {\em torsional rigidity} of the domain $\,\Omega^*\,$ is the value
$\,{\mathcal E} (\Omega^*) = \max_{u \in H_{1,c}^2 (\Omega^*)} \left( E_\Omega\, (u)\right)$. Since by (i) of Theorem \ref{symetrisation} 
$\,(f_\Omega)^* \in H_{1,c}^2 (\Omega^*, g^*) \,$ it follows
$${\mathcal E} (\Omega^*) \ge  E_{\Omega^*} \left( (f_\Omega)^*\right) \ge {\mathcal E} (\Omega).$$

\noindent Let us suppose that $\,{\mathcal E} (\Omega^*) = {\mathcal E} (\Omega)\,$, then
all the inequalities are equalities, in particular 
$$  \int_\Omega \,|\nabla f_\Omega|^2 \ dv_g  =  \alpha (M , M^*) \int_{\Omega^*} \,|\nabla (f_\Omega)^*|^2 \ dv_{g^*} $$
and $\,E_{\Omega^*} \left((f_\Omega)^*\right) = {\mathcal E} (\Omega^*) \,$.  Thus, since the set $\{f_\Omega>0\}$ coincides with the interior of $\Omega$, from the equality case of Theorem \ref{symetrisation} it follows that $\Omega^*\,$ is isometric to $\Omega$.

\vskip 0.3cm

\noindent {\it Proof of Theorem \ref{comparaison1}}:
\noindent As we have already pointed out in the introduction  the proof   follows immediately by Theorem \ref{comparaison}.
Indeed, Definition \ref{isoperimetre}  implies that $\,(M , g, x_0)$ is a 
PIMS for $(M,g)$ itself in the sense of the Definition \ref{espacemodele},
and  the constant $\,\alpha (M, M^*)\,$ is, in this case, always equal to $1$, 
because either $(M, g)$ 
has infinite volume, either the quotient of the volumes of the manifold and of the model space
is equal to $1$, because these two spaces coincide. $\square$
 \begin{remarque}\rm\label{explanation}
Let us return to the definition of the torsional rigidity
${\mathcal E} (\Omega) = E_\Omega  (f_\Omega)$. Two possible
definitions of this functional can be found in the classical literature: the one we considered here, i.e.
$${\mathcal E} (\Omega) =  \dfrac{1}{{\rm Vol} (\Omega)} \left(2 \int_\Omega  f_{\Omega} dv_g -\int_\Omega |\nabla f_\Omega|^2  dv_g\right)
=  \dfrac{1}{{\rm Vol} (\Omega)}  \int_\Omega f_\Omega dv_g $$
and, more frequently, the functional
$$\widetilde {\mathcal E} (\Omega) = 2 \int_\Omega \, f_\Omega\ dv_g
- \int_\Omega \,|\nabla f_\Omega|^2 \ dv_g =  \int_\Omega \, f_\Omega\ dv_g.$$
The critical or maximal domains (among all domains of prescribed volume) for the two functionals $\Omega \rightarrow {\mathcal E} (\Omega)$ and $\Omega \rightarrow \,\widetilde {\mathcal E} (\Omega)$ being the same, what is the interest of considering the first functional instead of the second one?
In fact, as $\,f_\Omega (x) \,$ is the 
\lq \lq mean exit time" for the paths of the Brownian motion issued from $\,x\,$,
$\,{\mathcal E} (\Omega) \,$ is the mean value of this \lq \lq mean exit time" with respect to all
possible initial points $\, x \in \Omega\,$, thus it still has some physical and stochastic 
meaning. 
Moreover if, on the same domain $\Omega$, we change the Riemannian metric $g$ in the homothetic metric $\lambda^2 g$, a direct computation gives:

\begin{equation}\label{homogeneite}
{\mathcal E} (\Omega , \lambda^2 \, g) = \lambda^2 \, {\mathcal E} (\Omega , g)\ \ \ \ \ {\rm and}
\ \ \ \ \ \widetilde {\mathcal E} (\Omega , \lambda^2 \, g) = \lambda^{n+2} \,
\widetilde  {\mathcal E} (\Omega , g),
\end{equation}
thus $\,{\mathcal E} (\Omega , g)\,$ has the same homogeneity as the Riemannian metric $\,g\,$.

But the main reason to prefer $\,{\mathcal E} (\Omega) \,$ to
$\,\widetilde {\mathcal E} (\Omega) \,$  is  Theorem
\ref{comparaison}, which provides a direct and simple comparison of the type 
$\,{\mathcal E} (\Omega) \le {\mathcal E} (\Omega^*) \,$, 
while, with respect to $\widetilde {\mathcal E}$, the comparison writes $\widetilde {\mathcal E}(\Omega) \leq \alpha(M, M^*) \widetilde {\mathcal E}(\Omega^*)$, 
where $\,\Omega\,$ is any domain on a Riemannian manifold (resp. on a compact Riemannian manifold $\,(M,g)\,)$ and 
$\,\Omega^*\,$ is a geodesic ball of the same volume (resp. of the same relative volume) 
on a model-space $\,(M^*, g^*, x^*)$.
\end{remarque}

\section{The proofs of Theorem \ref{compsphere2} and Theorem \ref{cheeger}}
In order to prove Theorem \ref{compsphere2} (see the end of this section)  which represents our main comparison between   torsional rigidities  in two different  compact manifolds (we remind the reader that the noncompact case has been treated in Corollary \ref{comparaisonCroke} of the Introduction), 
we need some known   isoperimetric inequalities (see Theorem \ref{Gromov}, Theorem  \ref{B-C-G} and  Theorem \ref{Perelman}) and, 
along the way, we also deduce some comparison results for the torsional rigidity (Corollary \ref{compsphere1}
and Corollary \ref{compsphere3} respectively).

\vskip 0.3cm

Revisiting Paul L\'{e}vy's work \cite{Le} (applied to convex bodies in the Euclidean space), M. Gromov (\cite{Gr}) 
proved the following celebrated isoperimetric inequality: 
\begin{theoreme}\label{Gromov} For every Riemannian manifold $\,(M,g)\,$ whose Ricci
curvature satisfies $\,{\rm Ric}_g \ge (n-1).g\,$, for every compact domain with smooth boundary 
$\,\Omega \,$ in $\, M \,$, let $\,\Omega^*\,$ be a geodesic ball of  the canonical sphere 
$\, (\mathbb S^n , g_0 )\,$ such that 
$\,\dfrac{{\rm Vol} (\Omega^*, g_0)}{{\rm Vol} (\mathbb S^n , g_0 )} 
= \dfrac{{\rm Vol} (\Omega , g )}{{\rm Vol} (M , g )}\,$, then
$$\dfrac{{\rm Vol}_{n-1}(\partial \Omega)}{{\rm Vol} (M , g )} \ge 
\dfrac{{\rm Vol}_{n-1}(\partial \Omega^*)}{{\rm Vol}  (\mathbb S^n , g_0 )}.$$
Moreover, this last inequality is an equality if and only if $\,\Omega\,$ is isometric to $\Omega^*$.
In other words, for any $\, x_0 \in \mathbb S^n \,$, $ \,(\mathbb S^n , g_0 , x_0)\,$ is a strict PIMS   for all the 
Riemannian manifolds $(M,g)$ which satisfy ${\rm Ric}_g \ge (n-1).g$.
\end{theoreme}

\begin{remarque}\rm
Theorem \ref{Gromov} is evidently sharp, because the canonical sphere $(\mathbb S^n , g_0 )$ satisfies its assumption 
\lq \lq $\,{\rm Ric}_g \ge (n-1).g\,$" (actually ${\rm Ric}_{g_0} = (n-1).g_0$), thus the Theorem applies when $\, (M, g) = 
(\mathbb S^n , g_0 )\,$, and because the isoperimetric inequality given by the Theorem
\ref{Gromov} is an equality when $\, (M, g) = (\mathbb S^n , g_0 )$ and when $\Omega$ is a geodesic ball of $(\mathbb S^n , g_0 )$.
\end{remarque}

Applying Theorems \ref{comparaison} and \ref{Gromov}, we obtain:

\begin{corollaire}\label{compsphere1}
For every Riemannian manifold $\,(M,g)\,$ whose Ricci
curvature satisfies $\,{\rm Ric}_g \ge (n-1).g\,$, for every compact domain with smooth boundary 
$\,\Omega \,$ in $\, M \,$, let $\,\Omega^*\,$ be a geodesic ball of the canonical sphere 
$\, (\mathbb S^n , g_0 )\,$ such that 
${{\rm Vol} (\Omega^*, g_0)}/{{\rm Vol} (\mathbb S^n , g_0 )} 
= {{\rm Vol}  (\Omega , g )}/{{\rm Vol} (M , g )}$, then
${\mathcal E} (\Omega) \leq {\mathcal E} (\Omega^*).$
Moreover, the equality $\,{\mathcal E} (\Omega) = {\mathcal E} (\Omega^*)\,$ is realized 
if and only if $\,\Omega\,$ is isometric to $\,\Omega^*\,$.
\end{corollaire}

\begin{remarque}\rm
It is easy to extend 
Corollary \ref{compsphere1} 
to every Riemannian manifold $\,(M,g)\,$ whose Ricci curvature satisfies $\,{\rm Ric}_g \ge 
K\,(n-1).g\,$ (with $\,K > 0\,$): in fact the Riemannian manifold $\,(M, K.g)\,$ then satisfies
$\,{\rm Ric}_{K.g} \ge (n-1). (K.g)\,$ and we can thus apply Theorem \ref{compsphere1} to the 
Riemannian manifold $\,(M, K.g)\,$; for every compact domain with smooth boundary 
$\,\Omega \,$ in $\, M \,$, if $\,\Omega^*\,$ is a geodesic ball on the Euclidean sphere 
$\, \mathbb S^n (\frac{1}{\sqrt{K}})\,$ of radius 
$\,\frac{1}{\sqrt{K}}\,$ and if $\,\Omega^{**}\,$ is a geodesic ball of  the canonical sphere 
$\, \mathbb S^n (1) = (\mathbb S^n , g_0 )\,$ such that 
$$\dfrac{{\rm Vol} (\Omega^{**})}{{\rm Vol} (\mathbb S^n , g_0 )} =
\dfrac{{\rm Vol} (\Omega^* )}{{\rm Vol} (\mathbb S^n (\frac{1}{\sqrt{K}}))} 
= \dfrac{{\rm Vol} (\Omega , g )}{{\rm Vol} (M , g )},$$ 
then, by (\ref{homogeneite}) and Theorem \ref{compsphere1},
$${\mathcal E} (\Omega , g) = \frac{1}{K}\  {\mathcal E} (\Omega , Kg) \le  \frac{1}{K}\  
{\mathcal E} (\Omega^{**} , g_0) = {\mathcal E} (\Omega^{**} , \frac{1}{K} g_0) = {\mathcal E} (\Omega^*),$$
where the last equality deduces from the fact that $\, (\mathbb S^n , \frac{1}{K} g_0 )\,$ is isometric to
$\mathbb S^n (\frac{1}{\sqrt{K}})$ and that this isometry maps $\Omega^{**}$ onto
$\Omega^{*}$. 
\end{remarque}

\begin{remarque}\rm
Corollary \ref{compsphere1}  is sharp
in the following sense:
for every $\,\beta \in ]0 , 1 [\,$, if we consider the set 
$\,{\mathcal W}_\beta\,$ of all domains $\,\Omega \,$ in all
the Riemannian manifolds $\,(M,g)\,$ whose Ricci curvature satisfies $\,{\rm Ric}_g \ge (n-1).g\,$
such that ${\rm Vol} (\Omega , g ) / {\rm Vol} (M , g ) = \beta\,$ then the geodesic 
ball $\,\Omega^*\,$ of the canonical sphere $\, (\mathbb S^n , g_0 )\,$ such that 
${\rm Vol} (\Omega^*, g_0)/ {\rm Vol} (\mathbb S^n , g_0 ) = \beta\,$ is an element
of $\,{\mathcal W}_\beta\,$.\\ 
If we consider the functional $\, \Omega \mapsto {\mathcal E} (\Omega) \,$
restricted to $\,{\mathcal W}_\beta\,$, then $\,\Omega^*\,$ is the point where this functional
attains its maximum.
\end{remarque}

Theorem \ref{Gromov} was improved and generalized to the case where the Ricci curvature has any sign
by P. B\'{e}rard, G. Besson and S. Gallot (\cite{B-B-G} Theorem (2) and \cite{Ga2}
Theorem 6.16 for a quantitatively improved version), who proved that

\begin{theoreme}\label{B-C-G} For any $\,K \in \mathbb R\,$, a PIMS  for all the 
$n$-dimensional Riemannian manifolds $(M,g)$ which satisfy ${\rm Ric}_g \ge (n-1)\,K .g$ 
and ${\rm diameter} (M,g) \leq D$ is given
by the Euclidean sphere of radius $\,R(K , D)\,$ (PIMS at any point) where $R(K, D)$ is defined by
 {\small $$ R(K , D) =  \left\{\begin{array}{l}
\frac{1}{\sqrt{K}} \left( \dfrac{\int_0^{\frac{D\sqrt{K}}{2}}\, (\cos t)^{n-1}\, dt}
{\int_0^{\frac{\pi}{2}}\, (\cos t)^{n-1}\, dt}\right)^{\frac{1}{n}} \ \ \ \ \ \ \ \ \ \ \  \ \ \ \ \ \ \  
{\rm if \ } K > 0\ \ \\ \\ 
\frac{n}{2} \left( \int_0^{\frac{\pi}{2}}\, (\cos t)^{n-1}\, dt \right)^{-\frac{1}{n}}\, D
 \ \ \ \ \ \ \ \ \ \ \ \ \ \ \  \ \ \ \ \ \ \  {\rm if \ } K = 0\ \ \\ \\
\frac{1}{\sqrt{|K|}} {\rm Max} \left(  \dfrac{\int_0^{D \sqrt{|K|}}\, 
(\cosh 2t)^{\frac{n-1}{2}}\, dt}{\int_0^{\pi}\, (\sin t)^{n-1}\, dt}\ ,\ 
\left( \dfrac{\int_0^{D \sqrt{|K|}}\, 
(\cosh 2t)^{\frac{n-1}{2}}\, dt}{\int_0^{\pi}\, (\sin t)^{n-1}\, dt}\right)^{\frac{1}{n}} \right)\\ 
\ \  {\rm if \ } K < 0
\end{array}\right.$$}
In other terms, for every compact domain with smooth boundary $\,\Omega \,$ in 
$\, M \,$, if $\,\Omega^*\,$ is a geodesic ball on the Euclidean sphere 
$\, \mathbb S^n (R(K , D))\,$ of radius 
$\,R(K , D)\,$ and if $\,\Omega^{**}\,$ is a geodesic ball of  the canonical sphere 
$\, \mathbb S^n (1) = (\mathbb S^n , g_0 )\,$ such that 
$$\dfrac{{\rm Vol} (\Omega^{**})}{{\rm Vol} (\mathbb S^n , g_0 )} =
\dfrac{{\rm Vol} (\Omega^* )}{{\rm Vol} (\mathbb S^n (R(K , D)))} 
= \dfrac{{\rm Vol} (\Omega , g )}{{\rm Vol} (M , g )}\ ,$$ 
then
\begin{equation}\label{sharpness}
\dfrac{{\rm Vol}_{n-1}(\partial \Omega)}{{\rm Vol} (M , g )} \ge 
\dfrac{{\rm Vol}_{n-1}(\partial \Omega^*)}{{\rm Vol} (\mathbb S^n (R(K, D)) )} =
\frac{1}{R(K , D)} \dfrac{{\rm Vol}_{n-1}(\partial \Omega^{**})}{{\rm Vol}  (\mathbb S^n , g_0 )}  .
\end{equation}
\end{theoreme}

\begin{remarque}\rm
This theorem is sharp, because the canonical sphere satisfies its assumptions 
\lq \lq $\,{\rm Ric}_g \ge (n-1).g\,$" and \lq \lq $\, {\rm diameter} \le \pi \,$",
 thus the Theorem \ref{B-C-G} applies to any domain $\Omega \subset \mathbb S^n$ when $\, (M, g) = (\mathbb S^n , g_0 )\,$, with the values
$\,K = 1\,$ and $\, D = \pi\,$ of the constants, and then the isoperimetric inequality (\ref{sharpness}) 
given by the theorem \ref{B-C-G} is an equality when $\, (M, g) = (\mathbb S^n , g_0 )\,$ and when $\Omega$ is a geodesic ball because, in this case, $\, R(K , D) = R(1, \pi) = 1\,$.
Moreover, under the assumptions \lq \lq $\,{\rm Ric}_g \ge (n-1).g\,$" and 
\lq \lq $\, (M, g) \,$ not isometric to $\,(\mathbb S^n , g_0 )\,$", Myers' theorem (and its equality case)
implies that $\,{\rm diameter} (M,g) < \pi\,$, and thus we can apply the Theorem \ref{B-C-G} 
with the values $\,K = 1\,$ and $\, D < \pi\,$ of the constants, which implies that, under these assumptions,
$\,R(K,D) < 1\,$. The isoperimetric inequality (\ref{sharpness}) is then strictly better than the one
of the canonical sphere.
Let us also remark that the smaller $\,R(K , D)\,$ is, the better is the isoperimetric 
inequality (\ref{sharpness}) given by the Proposition \ref{B-C-G}.
\end{remarque}

\begin{corollaire}\label{compsphere3}
For any $\,K \in \mathbb R\,$, for any $n$-dimensional Riemannian manifold $\,(M,g)\,$ 
which satisfies $\,{\rm Ric}_g \ge (n-1)\,K .g\,$ and $\, {\rm diameter} (M,g) \le D\,$, 
for every compact domain with smooth boundary $\,\Omega \,$ in $\, M \,$, if 
$\,\Omega^*\,$ is a geodesic ball on the Euclidean sphere 
$\, \mathbb S^n (R(K , D))\,$
and if $\,\Omega^{**}\,$ is a geodesic ball of  the canonical sphere 
$\, \mathbb S^n (1) = (\mathbb S^n , g_0 )\,$ such that 
$$\dfrac{{\rm Vol} (\Omega^{**})}{{\rm Vol} (\mathbb S^n , g_0 )} =
\dfrac{{\rm Vol} (\Omega^* )}{{\rm Vol} (\mathbb S^n (R(K , D)))} 
= \dfrac{{\rm Vol} (\Omega , g )}{{\rm Vol} (M , g )}\ ,$$ 
then
\begin{equation}
{\mathcal E} (\Omega) \le {\mathcal E} (\Omega^*) = R(K,D)^2 \ {\mathcal E} (\Omega^{**}).
\end{equation}
\end{corollaire}
\noindent
 {\it Proof}: Theorem \ref{B-C-G} shows that 
the Euclidean sphere $\, \mathbb S^n (R(K , D))\,$ of radius 
$R(K, D)$ is a PIMS for the Riemannian manifold $(M,g)$. 
For every compact domain with smooth boundary $\,\Omega \,$ in 
$\, M \,$, if $\,\Omega^*\,$ is a geodesic ball on the Euclidean sphere 
$\, \mathbb S^n (R(K , D))\,$ of radius $\,R(K , D)\,$ and if $\,\Omega^{**}\,$ 
is a geodesic ball of the canonical sphere 
$\, \mathbb S^n (1) = (\mathbb S^n , g_0 )\,$ such that 
$$\dfrac{{\rm Vol} (\Omega^{**})}{{\rm Vol} (\mathbb S^n , g_0 )} =
\dfrac{{\rm Vol} (\Omega^* )}{{\rm Vol} (\mathbb S^n (R(K , D)))} 
= \dfrac{{\rm Vol} (\Omega , g )}{{\rm Vol} (M , g )}\ ,$$ 
then Theorem \ref{comparaison} implies that 
$$
{\mathcal E} (\Omega) \le {\mathcal E} (\Omega^*) = R(K,D)^2 \,{\mathcal E} (\Omega^{**}),
$$
where the last equality deduces from the fact that the sphere of radius $\,R(K, D) \,$ is
isometric to $\,(\mathbb S^n , R(K,D)^2 . g_0 )\,$ and from
 formula (\ref{homogeneite}). $\ \square$

\vskip 0.5cm

We now recall an inequality due to G. Perelman \cite{Pe} (which is an improvement of a previous result of S. Ilias \cite{Il}).

\begin{theoreme}\label{Perelman}
Let $(M, g)$ be a $n$-dimensional  compact Riemannian manifold.
Assume that $M$ is {\em not} diffeomorphic to $\mathbb S^n$, that ${\rm Ric}_g \ge (n-1).g$ and that the sectional curvature of $(M , g)$ is  $ \geq - \kappa^2$.
Then there exists a constant  $\varepsilon (n , \kappa) > 0$
such that ${\rm diameter} (M,g) \leq \pi - \varepsilon (n , \kappa)$.
\end{theoreme}

\begin{remarque}\rm
By applying  the Theorem \ref{B-C-G} with the values $\,K = 1\,$ and $\, D = \pi - 
\varepsilon (n , \kappa)\,$ of the constants, which implies that, under these assumptions,
\begin{equation}\label{gap}
R(K,D)  = R(1, \pi - \varepsilon (n , \kappa) ) =  
\left( 1 - \dfrac{\int_0^{\frac{\varepsilon (n , \kappa)}{2}}\, (\sin t)^{n-1}\, dt}
{\int_0^{\frac{\pi}{2}}\, (\sin t)^{n-1}\, dt}\right)^{\frac{1}{n}}, 
\end{equation}
we observe that,  with respect to the isoperimetric inequality of the canonical sphere, the isoperimetric inequality
on $\, (M, g)\,$ induced by (\ref{sharpness}) is improved by some factor which is bounded far from $1$.
\end{remarque}

\vskip 0.3cm

We are now in the position to prove Theorem  \ref{compsphere2}
(notice that it improves Corollary \ref{compsphere1}).

\vskip 0.3cm

\noindent {\it Proof of Theorem \ref{compsphere2}}\,:  Applying Theorem \ref{B-C-G} (with the values $\,K = 1\,$ and $\, D = 
{\rm diameter} (M,g) \,$ of the constants) we prove that the Euclidean sphere 
$\, \mathbb S^n (R(1 , D))\,$ of radius 
$\,R(1, D) = R(1,{\rm diameter} (M,g))\,$ is a  PIMS (at any point) for the Riemannian manifold $(M,g)$. 
For every compact domain with smooth boundary $\,\Omega \,$ in 
$\, M \,$, if $\,\Omega^0\,$ is a geodesic ball on the Euclidean sphere 
$\, \mathbb S^n (R(1 , D))\,$ of radius $\,R(1 , D)\,$ and if $\,\Omega^*\,$ 
is a geodesic ball of  the canonical sphere 
$\, \mathbb S^n (1) = (\mathbb S^n , g_0 )\,$ such that 
$$\dfrac{{\rm Vol} (\Omega^*)}{{\rm Vol} (\mathbb S^n , g_0 )} =
\dfrac{{\rm Vol} (\Omega^0 )}{{\rm Vol} (\mathbb S^n (R(1 , D)))} 
= \dfrac{{\rm Vol} (\Omega , g )}{{\rm Vol} (M , g )}\ ,$$ 
then  Theorem \ref{comparaison} implies that 
\begin{equation}\label{gap1}
{\mathcal E} (\Omega) \le {\mathcal E} (\Omega^0) = R(1,D)^2 \,{\mathcal E} (\Omega^*),
\end{equation}
where the last equality deduces from the fact that the sphere of radius $\,R(1, D) \,$ is
isometric to $\,(\mathbb S^n , R(1,D)^2 . g_0 )\,$ and from
 formula (\ref{homogeneite}).

Let us first suppose that $(M, g)$ is not isometric to $(\mathbb S^n , g_0 )$,
then Myers' theorem (and its equality case) implies that $\,{\rm diameter} (M,g) < \pi\,$, 
and thus that $\,R(1, D) < 1\,$ (if $\, D = {\rm diameter} (M,g) \,$)
by the definition of $\,R(K, D)\,$. Using the fact that $\,R(1, D) < 1\,$
in the inequality (\ref{gap1}), we conclude that, if $\, (M, g) \,$ is not isometric to 
$\,(\mathbb S^n , g_0 )\,$, then $\,{\mathcal E} (\Omega) < {\mathcal E} (\Omega^*)\,$ for every 
compact domain with smooth boundary $\,\Omega \,$ in $\, M \,$, which proves the
part (i) of the Theorem \ref{compsphere2}.

If we now suppose that $M$ is not diffeomorphic to $\,\mathbb S^n\,$, we
know, by Theorem \ref{Perelman}, that $\,{\rm diameter} (M,g) \leq \pi - \varepsilon (n , \kappa)\,$ in this case
and thus that we can choose the value 
$D = \pi - \varepsilon (n , \kappa)\,$ for the upper bound of the diameter of $(M, g)$.
Using the inequality (\ref{gap1}) (and the formula (\ref{gap}) for the computation of
$R(1, \pi - \varepsilon (n , \kappa) )$, we get that
$${\mathcal E} (\Omega) \leq  \left( 1 - \dfrac{\int_0^{\frac{\varepsilon (n , \kappa)}{2}} (\sin t)^{n-1} dt}
{\int_0^{\frac{\pi}{2}}\, (\sin t)^{n-1}\, dt}\right)^{\frac{2}{n}}{\mathcal E} (\Omega^*)$$
for every compact domain with smooth boundary $\,\Omega \,$ in $M$, which proves 
part (ii) of  Theorem \ref{compsphere2}.  $\square$

In order to prove Theorem \ref{cheeger} we need the following:

\begin{lemme}\label{lemmacheeger}
Let $(M, g)$ be a Riemannian manifold. Then, for any compact domain $\Omega\subset M$ and for any smooth nonnegative function 
$f$ on $\Omega$ which vanishes on $\partial\Omega$, one has:
$$\int_{\Omega}fdv_g=\int_0^{\sup f}A(t)dt,$$
where $A(t)=\Vol (\Omega_t)$ and  where $\Omega_t$ denotes the set of points $x\in\Omega$ such that $f(x)>t$.
\end{lemme}
\noindent {\it Proof}: 
Let $t_i=\frac{i}{N}\sup f$ (for every $i\in \{0, \dots N\}$).
The function $\, t \mapsto A(t) = {\rm Vol} (\Omega_t)\,$ being 
strictly decreasing, we have
\begin{equation}\label{tAt}
\sum_{i = 0}^{N-1} t_i  \left( A(t_i) - A(t_{i+1}) \right) \le \int_\Omega \, f\, dv_g 
\le \sum_{i = 0}^{N-1} t_{i+1}  \left( A(t_i) - A(t_{i+1}) \right) .
\end{equation}
Let   $S^+_N$ (resp. $S^{-}_N$) denote the right (resp. the left)
hand side of (\ref{tAt}).
This is an approximation from above (resp. from below) of the integral 
$\int_0^{\sup f}A(t)dt$.
As $0\leq S^+_N-S^-_N\leq\frac{\sup f}{N}A(0)$, when $N\rightarrow\infty$, 
$S^+_N-S^-_N\rightarrow 0_+$ and $S^+_N$, $S^-_N$ both go to 
$\int_0^{\sup f}A(t)dt$ (and to $\int_{\Omega}fdv_g$ by (\ref{tAt})).
$\square$

\vskip 0.3cm

\noindent {\it Proof of Theorem \ref{cheeger}}:
By the definition of ${\mathcal E}(\Omega)$ and by Lemma \ref{lemmacheeger}, we have
\begin{equation}\label{torsvol}
{\rm Vol}(\Omega) \ {\mathcal E}(\Omega)=\int_\Omega f_{\Omega}dv_g=\int_{[0, \sup f_\Omega]\setminus S(f)}A(t)dt,
\end{equation}
where $A(t)={\rm Vol} (\Omega_t)$ and where $\Omega_t$  denotes the set of points $x\in\Omega$ such that $f_\Omega(x)>t$.
For every regular value  $t$ of $f_{\Omega}$ one has:
$$A(t)\leq {\rm Vol}(\Omega)\leq {\rm Vol}(M, g)/2$$
and thus, by the definition of Cheeger's isoperimetric constant,
$${\rm Vol}_{n-1}(\partial\Omega_t)\geq H(M, g)A(t).$$
From this and from (\ref{torsvol}) we deduce
$${\rm Vol}(\Omega) \ {\mathcal E}(\Omega)\leq \frac{1}{H(M, g)}\int_{[0, \sup f_\Omega]\setminus S(f)}{\rm Vol}_{n-1}(\partial\Omega_t)dt=
\frac{1}{H(M, g)}\int_\Omega|\nabla f_\Omega |dv_g,$$
where, in the last equality, we have used the coarea formula (\ref{coarea}).
Thus, by Cauchy--Schwarz inequality
$${\rm Vol}(\Omega) \ {\mathcal E}(\Omega)\leq\frac{1}{H(M, g)}\left({\rm Vol}(\Omega)\right)^{\frac{1}{2}}\left(\int_\Omega|\nabla f_\Omega |^2dv_g\right)^{\frac{1}{2}}$$
and hence, since  ${\mathcal E}(\Omega)=\frac{1}{{\rm Vol(\Omega)}}\int_\Omega|\nabla f_\Omega |^2dv_g$, one gets 
$\left({\mathcal E}(\Omega)\right)^{\frac{1}{2}}\leq\frac{1}{H(M, g)}$ which ends the proof of the theorem.
$\square$

\subsection* {On the sharpness of Theorem \ref{cheeger}}
The following  proposition  shows that Theorem \ref{cheeger} is sharp.
\begin{proposition}\label{propsharp}
There exists  a family of  $n$-dimensional compact Riemannian manifolds $(M, g_{\varepsilon})$, $0<\varepsilon<1$,
such that for every $\beta \in [\frac{\varepsilon}{4}, \frac{1}{2}]$
there exists a compact domain $\Omega \subset M$ 
and a universal constant $B=B(n)$ such that
\begin{equation}\label{volbeta}
\frac{{\rm Vol}(\Omega , g_{\varepsilon})}{{\rm Vol}(M, g_{\varepsilon})}=\beta
\end{equation}
and
\begin{equation}\label{eomegahm}
{\mathcal E}(\Omega, g_{\varepsilon})\geq \frac{(1-\varepsilon)(1-B\sqrt{\varepsilon})^2}{H(M, g_{\varepsilon})^2} .
\end{equation}
\end{proposition}
\noindent {\it Proof}: 
Let  $R$,  $\delta$ and $\varepsilon$ be  positive real numbers 
such that $\delta >\frac{1}{n}$ and $\varepsilon \leq \frac{1}{n\delta}<1$.
 Define two other  positive real numbers $\lambda$ and $\eta$ by: 
$$\lambda =\frac{e^{\delta R}}{\varepsilon}\sqrt{1-\varepsilon^2\delta^2e^{-2\delta R}},\ \ \ \  \eta =\frac{1}{\lambda}\arctan (\frac{\lambda}{\delta}).$$
Consider the compact  manifold $M$ obtained as the quotient of  $[-(R+\eta), R+\eta]\times \mathbb S^{n-1}$  by   identifying 
all  the points $\{R+\eta\}\times \mathbb S^{n-1}$ (resp. $\{-(R+\eta)\}\times \mathbb S^{n-1}$) to a single point $x_0$ (resp. $x_1$)
(see Subsection \ref{Example} above). 
Denote by  
$$\pi:  [-(R+\eta), R+\eta]\times \mathbb S^{n-1}\rightarrow M$$
the  corresponding quotient map.
We endow $M$ with the metric $g_\varepsilon$ defined at the point $(t, v)\in [-(R+\eta), R+\eta]\times \mathbb S^{n-1}$ by:
$$g_{\varepsilon}=(dt)^2+b_{\varepsilon}(t)^2g_0,$$
where  $b_{\varepsilon}$ is the even function  on $[-(R+\eta), R+\eta]$ given  by
$$b_{\varepsilon}(t)= \left\{\begin{array}{l}
\varepsilon  e^{-\delta t}\ \  \ \ \ \ \ \ \ \ \ \ \ \ \ \ \ \ \ \ \ \ {\rm if \ }\ \ \ t\in [0, R], \\
\varepsilon e^{-\delta R} \left(\frac{\sin[\lambda (\eta +R-t)]}{\sin (\lambda\eta)}\right) \ \ \ {\rm if \ }\ \ \ t\in [R, R+\eta].
\end{array}\right.$$
It is easily seen that $b_{\varepsilon}$ is $C^{1}$ on $]0, R+\eta [$ and  that 
$$b_{\varepsilon}(-(R+\eta))=b_{\varepsilon}(R+\eta)=0, \  b^{'}(R+\eta)=-1,  b^{'}(-(R+\eta))=1.$$
Moreover we have
\begin{equation}\label{bepsminore}
b_{\varepsilon}(t)\leq \varepsilon e^{-\delta t}, \ t\in [0, R+\eta].
\end{equation}

For $r\in [0, R+\eta]$  consider the compact domain
$$\Omega_r=\pi ( [r, R+\eta]\times \mathbb S^{n-1})\subset M,$$
we claim that   
\begin{equation}\label{ineqvolvoldelta}
\frac{{\rm Vol}_{n-1}(\partial\Omega_r, g_{\varepsilon})}{{\rm Vol} (\Omega_r, g_{\varepsilon})}\geq (n-1)\delta, \ r\in [0, R+\eta].
\end{equation}
Indeed, on the one hand,   when  $r\in [0, R]$, 
(\ref{bepsminore}) 
yields 
$$\frac{{\rm Vol}_{n-1}(\partial\Omega_r, g_{\varepsilon})}{{\rm Vol} (\Omega_r, g_{\varepsilon})}=\frac{b_{\varepsilon}(r)^{n-1}}{\int_r^{R+\eta}b_{\varepsilon}(t)^{n-1}dt}
\geq \frac{e^{-(n-1)\delta r}}{\int_r^{R+\eta}e^{-(n-1)\delta t}dt}
\geq (n-1)\delta .$$
On the other  hand,  when $\ r\in [R, R+\eta]$, by setting $\tilde r=R+\eta -r$ we get
$$\frac{{\rm Vol}_{n-1}(\partial\Omega_r, g_{\varepsilon})}{{\rm Vol} (\Omega_r, g_{\varepsilon})}=\frac{\sin[\lambda (R+\eta -r)]^{n-1}}{\int_r^{R+\eta}\sin[\lambda (R+\eta -t)]^{n-1}dt}
\geq \frac{\lambda\sin[\lambda \tilde r]^{n-1}}{\int_0^{\lambda\tilde r}(\sin t)^{n-1}dt}$$
$$\geq \frac{\lambda}{\int_0^{\frac{\pi}{2}}(\sin t)^{n-1}dt}\geq \lambda
\geq (n-1)\delta ,$$
where the second inequality comes from the fact that the function $\xi \mapsto 
\frac{(\sin \xi )^{n-1}}{\int_0^{\xi}(\sin t)^{n-1}dt}$
is decreasing on the interval $[0, \frac{\pi}{2}]$ 
and $\lambda \tilde r\leq \lambda \eta =\arctan (\frac{\lambda}{\delta})<\frac{\pi}{2}$
and where the last inequality deduces from the definition of $\lambda$
and from the assumption
$\varepsilon \leq \frac{1}{n\delta}$. Hence the claim  (\ref{ineqvolvoldelta}) is proved.
Therefore, if we define
$$H_{\rm rad}(M, g_{\varepsilon}):=\inf_{r\in [0, R+\eta[}\frac{{\rm Vol}_{n-1}(\partial \Omega_r, g_{\varepsilon})}{{\rm Vol}(\Omega_r, g_{\varepsilon})}$$
we get 
$$H_{\rm rad}(M, g_{\varepsilon})\geq  (n-1)\delta .$$
By the method developed in Appendix A.4 of \cite{Ga2}, we can find a universal constant $B=B(n)$
such that
\begin{equation}\label{16}
(1-B\sqrt{\varepsilon})(n-1)\delta\leq (1-B\sqrt{\varepsilon})H_{\rm rad} (M, g_{\varepsilon})\leq H(M, g_{\varepsilon}).
\end{equation}

Consider now the test function $f:\Omega_r\rightarrow \real^+$ 
$$f(t, v)=u(t)= \left\{\begin{array}{l}
\frac{t-r}{(n-1)\delta}\ \  \ {\rm if \ }\ \ \ t\in [r, R], \\
\frac{R-r}{(n-1)\delta} \ \ \ {\rm if \ }\ \ \ t\in [R, R+\eta].
\end{array}\right.$$
By a straightforward computation we get:
\begin{equation}\label{dieci}
\int_{\Omega_r}|\nabla f|^2dv_{g_\varepsilon}<\frac{1}{(n-1)^2\delta^2} {\rm Vol}(\Omega_r, g_\varepsilon), 
\end{equation}
and,  using  (\ref{bepsminore}), and setting $\omega_{n-1}={\rm Vol}(\mathbb S^{n-1}, g_0)$
$$\int_{\Omega_r} f^2dv_{g_{\varepsilon}}\leq \frac{\omega_{n-1}\varepsilon^{n-1}} {(n-1)^2\delta^2}\int_{r}^{R+\eta}(t-r)^2e^{-(n-1)\delta t}dt<+\infty .$$
As $f$ is piecewise $C^1$ and vanishes on $\partial\Omega_r$ (because $u(r)=0$), we get that 
\begin{equation}
f\in H_{1, c}^2(\Omega_r).
\end{equation}
Using  (\ref{bepsminore}) we also have:
\begin{equation}\label{bintegrating}
\frac{\omega_{n-1}\varepsilon^{n-1}}{(n-1)\delta}e^{-(n-1)\delta r}\left(1-e^{-(n-1)\delta(R-r)}\right)\leq {\rm Vol}(\Omega_r, g_\varepsilon)\leq 
\frac{\omega_{n-1}\varepsilon^{n-1}}{(n-1)\delta}e^{-(n-1)\delta r}
\end{equation}
Integrating by parts, we get:
$$\int_{\Omega_r}fdv_{g_{\varepsilon}}\geq \frac{\omega_{n-1}\varepsilon^{n-1}}{(n-1)\delta}\int_r^R(t-r)e^{-(n-1)\delta t}dt
\geq $$
$$\geq  \frac{\omega_{n-1}\varepsilon^{n-1}}{(n-1)^3\delta^3}e^{-(n-1)\delta r}\left[1-\left(1+(n-1)\delta (R-r)\right)e^{-(n-1)\delta (R-r)}\right]$$
Combining this last equality with (\ref{bintegrating}) it follows that:
\begin{equation}\label{aintegrating}
\frac{1}{{\rm Vol}(\Omega_r, g_\varepsilon)}\int_{\Omega_r}fdv_{g_{\varepsilon}}\geq \frac{1-\left(1+(n-1)\delta(R-r)\right)e^{-(n-1)\delta (R-r)}}{(n-1)^2\delta^2}.
\end{equation}
As $f\in H_{1, c}^2(\Omega_r)$ one may apply Definition  \ref{fonctionnelle} which gives (with the help of (\ref{aintegrating}) and (\ref{dieci})):
$${\mathcal E}(\Omega_r, g_{\varepsilon})\geq \frac{2\int_{\Omega_r}fdv_{g_\varepsilon}-
\int_{\Omega_r}|\nabla f|^2dv_{g_\varepsilon}}{{\rm Vol}(\Omega_r, g_\varepsilon)}
\geq \frac{1-2\left[1+(n-1)\delta(R-r)\right]e^{-(n-1)\delta (R-r)}}{(n-1)^2\delta^2}.$$
Let $A_{\varepsilon}$ be the unique solution of the equation 
$$2(1+x)e^{-x}=\varepsilon .$$
For every $R\geq \frac{2A_{\varepsilon}}{(n-1)\delta}$ and $0\leq r\leq \frac{R}{2}$ we have
$(n-1)\delta (R-r)\geq A_{\varepsilon}$ and thus
\begin{equation}\label{eq2eps}
2\left[1+(n-1)\delta(R-r)\right]e^{-(n-1)\delta (R-r)}\leq\varepsilon .
\end{equation}
Using (\ref{16}) we thus get (when $R\geq \frac{2A_{\varepsilon}}{(n-1)\delta}$ and $0\leq r\leq \frac{R}{2}$):
$${\mathcal E}(\Omega_r, g_{\varepsilon})\geq \frac{1-\varepsilon}{(n-1)^2\delta^2}\geq \frac{(1-\varepsilon )(1-B\sqrt{\varepsilon})^2}{H(M, g_{\varepsilon})^2},$$
which proves (\ref{eomegahm}).
As $x\mapsto \frac{x}{1+x}$ is increasing on $[0,+\infty [$, we have (again by (\ref{bepsminore})):
$$v(r):=\frac{{\rm Vol}(\Omega_r , g_{\varepsilon})}{{\rm Vol}(M, g_{\varepsilon})}=\frac{\int_r^{R+\eta}b_{\varepsilon}(t)^{n-1}dt}{2\int_0^{R+\eta}b_{\varepsilon}(t)^{n-1}dt}
\leq\frac{\int_r^{+\infty}e^{-(n-1)\delta t}dt}{2\int_0^{+\infty}e^{-(n-1)\delta t}dt}\leq \frac{1}{2} e^{-(n-1)\delta r}.$$
As $v(0)=\frac{1}{2}$ in order to prove (\ref{volbeta}), it is enough to show  that $v(\frac{R}{2})\leq\frac{1}{2}e^{-(n-1)\delta \frac{R}{2}}\leq\frac{\varepsilon}{4}$
which follows immediately from (\ref{eq2eps}) with $r=\frac{R}{2}$. This concludes the proof of the proposition.
$\ \square$

\end{document}